\input{psfig}
\documentstyle[11pt]{article}
\input amssym.def
\input amssym.tex
\newcommand{\beql}[1]{\begin{equation}\label{#1}}
\newcommand{\eeq}{\end{equation}}
\newcommand{\comment}[1]{}

\newcommand{\eqref}[1]{(\ref{#1})}

\newcommand{\Abs}[1]{{\left|{#1}\right|}}

\newcommand{\Norm}[1]{{\left\|{#1}\right\|}}

\newcommand{\Qed}{\newline\mbox{$\blacksquare$}}

\newcommand{\Set}[1]{{\left\{{#1}\right\}}}

\newcommand{\ToAppear}[1]{\raisebox{15mm}[10pt][0mm]{\makebox[0mm]{\makebox[\textwidth][r]{\small #1}}}}

\newcommand{\RR}{{\Bbb R}}
\newcommand{\CC}{{\Bbb C}}
\newcommand{\ZZ}{{\Bbb Z}}
\newcommand{\NN}{{\Bbb N}}

\newcounter{open}
\newcounter{mysec}
\def\themysec{\arabic{mysec}}
\newcommand{\mysection}[1]{
  \vskip 0.25in
  \refstepcounter{mysec}\centerline{\large\bf \S\themysec.\ {#1}}
  \addcontentsline{toc}{section}{{\bf \themysec.}\ {#1}}
}

\newcounter{mysubsec}[mysec]
\newtheorem{theorem}{Theorem}

\newtheorem{lemma}{Lemma}

\begin{document}

\begin{center}
{\Large \bf On the uniform distri\ToAppear{\tt PREPRINT}bution
in residue classes}\\
{\Large \bf of dense sets of integers with distinct sums.}\\
\ \\
{\sc Mihail N. Kolountzakis\footnote{
Department of Mathematics, University of Crete,
714 09 Iraklio, Greece. E-mail: {\tt kolount@math.uch.gr}}
\footnote{Partially supported by the U.S. National Science Foundation,
under grant DMS 97-05775.}}\\
\ \\
\small July 1998 %
\end{center}

\begin{abstract}
A set ${\cal A} \subseteq \Set{1,\ldots,N}$ is of type $B_2$ if all sums
$a+b$, with $a\ge b$, $a,b\in {\cal A}$, are distinct.
It is well known that the largest such set is of size asymptotic to $N^{1/2}$.
For a $B_2$ set ${\cal A}$ of this size we show that, under mild assumptions
on the size of the modulus $m$ and on the difference $N^{1/2}-\Abs{{\cal A}}$
(these quantities should not be too large) the elements of ${\cal A}$
are uniformly distributed in the residue classes mod $m$.
Quantitative estimates on how uniform the distribution is are also provided.
This generalizes recent results of Lindstr\"om
whose approach was combinatorial.
Our main tool is an upper bound on the minimum of a cosine sum of $k$ terms,
$\sum_1^k \cos{\lambda_j x}$,
all of whose positive integer frequencies $\lambda_j$
are at most $(2-\epsilon)k$ in size.
\end{abstract}

\mysection{Introduction and results}
\label{sec:intro}

\noindent
A set ${\cal A} \subseteq \Set{1,\ldots,N}$ is of type $B_2$ if all sums
$$
a+b,\ \ \mbox{with $a\ge b,\ a,b\in {\cal A}$},
$$
are distinct.
(Such sets are also called {\em Sidon}, but the term has a very different
meaning in harmonic analysis.)
This is easily seen to be equivalent to all differences
$a-b$, with $a\neq b,\ a,b\in {\cal A}$, being distinct.
It is an old theorem of Erd\H os and Tur\'an \cite{ET41,HR83,K96}
that the size of the largest $B_2$ subset of $\Set{1,\ldots,N}$
is at most $N^{1/2} + O(N^{1/4})$.
It is also known \cite{BC63,HR83} that there exist $B_2$ subsets
of $\Set{1,\ldots,N}$ of size $\sim N^{1/2}$.

In this note we consider such dense $B_2$ subsets ${\cal A}$ of $\Set{1,\ldots,N}$,
i.e., sets of size $N^{1/2} + o(N^{1/2})$, and prove,
under mild conditions on $\Abs{{\cal A}}$ and the modulus $m$, which is also
allowed to vary with $N$, that they are uniformly distributed mod $m$.
More precisely, let
$$
a(x) = a_m(x) = \Abs{\Set{a \in {\cal A}: a = x \bmod m}},
\ \ \mbox{for $x \in \ZZ_m$},
$$
be the number of elements of ${\cal A}$ with residue $x \bmod m$.
We shall show, for example,
that if $\Abs{{\cal A}} \sim N^{1/2}$ and $m$ is a constant then, as $N\to\infty$,
\beql{unif}
a(x) = {\Abs{{\cal A}}\over m} + o\left({\Abs{{\cal A}}\over m}\right).
\eeq
We shall also obtain bounds on the error term.
These bounds will depend on $\Abs{{\cal A}}, m$ and $N$.

Previously Lindstr\"om \cite{L98} showed precisely \eqref{unif}
using a combinatorial method, thus answering a question posed in \cite{ESS94}.
Under the additional assumptions
\beql{assumptions}
m=2 \ \ \mbox{and}\ \ \Abs{{\cal A}} \ge N^{1/2}
\eeq
he obtained the bound $O(N^{3/8})$ for the error term in \eqref{unif}.

Here we use an analytic method which has previously been used
\cite{K96} to prove and generalize the Erd\H os-Tur\'an theorem mentioned
above.
The core of our technique is the following theorem \cite{K96}
which was proved in connection with the so called {\em cosine problem}
of classical harmonic analysis.
\begin{theorem} \label{th:cos}
Suppose $0\le f(x) = M + \sum_1^N\cos{\lambda_j x}$, with
the integers $\lambda_j$ satisfying
$$
1\le\lambda_1<\cdots<\lambda_N \le (2-\epsilon)N,
$$
for some $\epsilon>3/N$.
Then
\begin{equation} \label{eq:main}
M > A \epsilon^2 N,
\end{equation}
for some absolute positive constant $A$.
\end{theorem}
Our main theorem, of which Lindstr\"om's result is a special case,
is the following.
\begin{theorem}\label{th:main}
Suppose ${\cal A} \subseteq \Set{1,\ldots,N}$ is a $B_2$ set and that
$$
k = \Abs{{\cal A}} \ge N^{1/2} - \ell(N),
\ \ \mbox{with $\ell(N) = o(N^{1/2})$}.
$$
Assume also that $m = o(N^{1/2})$.
Then we have
\beql{result}
\Norm{a(x) - {k\over m}}_2 \le 
 C \left\{ \begin{array}{cl}
{N^{3/8} \over m^{1/4}} & \mbox{if $\ell \le N^{1/4}m^{1/2}$}\\
\ & \ \\
{N^{1/4}\ell^{1/2}\over m^{1/2}} & \mbox{else}.
\end{array}\right.
\eeq
\end{theorem}
(In our notation $\ell$ need not be a positive quantity. If it is
negative (i.e., $k>N^{1/2}$) the first of the two alternatives
holds in the upper bound.)

We use the notation
$\Norm{f}_p = \left(\sum_{x\in\ZZ_m}\Abs{f(x)}^p\right)^{1/p}$, for
$f:\ZZ_m\to\CC$ and $1\le p <\infty$, and also
$\Norm{f}_\infty = \max_{x\in\ZZ_m}\Abs{f(x)}$.
We obviously have $\Norm{f}_\infty \le \Norm{f}_p$, for all $f$
and $1\le p < \infty$.

\noindent{\bf Remarks.}
It follows easily from Theorem \ref{th:main} that in the following
two cases we have uniform distribution in residue classes mod $m$.
\begin{enumerate}
\item
When $\ell \le N^{1/4}m^{1/2}$ and $m=o(N^{1/6})$ we have
\beql{unif-case-1}
\Norm{a(x)-{k\over m}}_\infty \le
\Norm{a(x)-{k\over m}}_2 \le C{N^{3/8}\over m^{1/4}} = o\left({k\over m}\right).
\eeq
\item
When $\ell \ge N^{1/4}m^{1/2}$
and $m=o\left({N^{1/2}\over \ell}\right)$ we have
\beql{unif-case-2}
\Norm{a(x)-{k\over m}}_\infty \le
\Norm{a(x)-{k\over m}}_2 \le C {N^{1/4}\ell^{1/2}\over m^{1/2}} =
 o\left({k\over m}\right).
\eeq
\end{enumerate}
In these two cases we have uniform distribution ``in the $\ell^2$ sense''
as well as in the $\ell^\infty$ sense.

As a comparison to the result that Lindstr\"om obtained under assumptions
\eqref{assumptions}, we obtain that whenever $m$ is a constant and
$\ell \le C N^{1/4}$ we have 
$$
\Norm{a(x)-{k\over m}}_2 \le C {N^{1/4}\ell^{1/2}\over m^{1/2}}
 \le C_m N^{3/8}.
$$
As is customary, $C$ denotes an absolute positive constant, not necessarily
the same in all its occurences, while $C$ subscripted is allowed to depend
only on the subscripts.

\mysection{Proofs}
\label{sec:proofs}

For the proof of Theorem \ref{th:main}
we shall need the following two lemmas, the first of which is
elementary and the second a consequence of Theorem \ref{th:cos}.
\begin{lemma}\label{lm:cst}
If $a:\ZZ_m\to\CC$ and $S = \sum_{x\in\ZZ_m}a(x)$ then
$$
\sum_{x\in\ZZ_m}\Abs{a(x) - {S\over m}}^2 =
 \sum_{x\in\ZZ_m}\Abs{a(x)}^2 - {S^2\over m}.
$$
\end{lemma}
{\bf Proof.}
Let $a(x) = {S\over m} + \delta(x)$ for $x\in\ZZ_m$.
It follows that $\sum_{x\in\ZZ_m}\delta(x) = 0$.
Then
\begin{eqnarray*}
\sum_{x\in\ZZ_m} \Abs{a(x)}^2 &=&
 \sum_{x\in\ZZ_m}\left( {S^2\over m^2} + \Abs{\delta(x)}^2 +
	2 {S\over m} {\rm Re\,}\delta(x) \right) \\
 &=& {S^2\over m} + \sum_{x\in\ZZ_m} \Abs{\delta(x)}^2.
\end{eqnarray*}
\Qed

\begin{lemma}\label{lm:cos_m}
Suppose $\lambda_j \in \NN, j=1,\ldots,N$, are distinct positive integers
and define
$$
N_m = \Abs{\Set{\lambda_j:\ \lambda_j = 0 \bmod m}}.
$$
If
$$
0 \le p(x) = M + \sum_{j=1}^N \cos{\lambda_j x},\ \ \ (x \in \RR),
$$
and
$$
\lambda_j \le (2-\epsilon) N_m m,\ \ \mbox{for all $\lambda_j=0 \bmod m$},
$$
for some $\epsilon> 3/N_m$,
then we have
$$
M > A \epsilon^2 N_m,
$$
for some absolute positive constant $A$.
\end{lemma}
{\bf Proof.}
The measure $\mu$ on $[0,2\pi)$ with $\widehat{\mu}(n)=1$ if $m$ divides $n$
and $\widehat{\mu}(n)=0$ otherwise is nonnegative.
Let
$$
q(x) = p(x)\star\mu = M + \sum_{m \mid \lambda_j} \cos{\lambda_j x} \ge 0.
$$
Define also the polynomial
$$
r(x) = q\left({x \over m}\right)
 =  M + \sum_{m \mid \lambda_j} \cos{{\lambda_j\over m}x}
 \ge 0.
$$
By Theorem \ref{th:cos} and the assumption
$$
{\lambda_j \over m} \le (2-\epsilon)N_m
$$
we get $M\ge A \epsilon^2 N_m$ as desired.
\Qed

\noindent
{\bf Proof of Theorem \ref{th:main}.}
Write
$$
d(j) = \Abs{\Set{(a,b) \in {\cal A}^2:\ a-b = j \bmod m}},\ \ \ (j\in\ZZ_m),
$$
and notice that, by the Cauchy-Schwarz inequality,
$$
d(j) = \sum_{i\in\ZZ_m} a(i)a(i+j) \le \sum_{i\in\ZZ_m} (a(i))^2 = d(0),
\ \ \ (j\in\ZZ_m).
$$
We also clearly have $\sum_{i\in\ZZ_m} d(i) = k^2$ which implies
$$
d(0) \ge {k^2 \over m}.
$$
Define the nonnegative polynomial
\begin{eqnarray*}
f(x) &=& \Abs{\sum_{a \in {\cal A}} e^{i a x}}^2\\
 &=& k + \sum_{a\neq b,\ a,b \in {\cal A}} e^{i(a-b)x}\\
 &=& k + 2 \sum_j \cos{\lambda_j x},
\end{eqnarray*}
where the set $\Set{\lambda_j}$ consists of all differences
$a-b$, with $a,b\in {\cal A}$, $a>b$, which are all distinct since ${\cal A}$ is
of type $B_2$.
(Notice that $1\le \lambda_j \le N$.)
With the notation of Lemma \ref{lm:cos_m} we have
$$
d(0) = k + 2 N_m.
$$
Since $k \sim N^{1/2}$ and $m = o(N^{1/2}) = o(k)$ we may suppose that,
for $N$ large enough,
$$
{1\over 2}N^{1/2} < k < 2 N ^{1/2}
$$
and
$$
m < {1\over2}k.
$$
Hence
$$
{3\over N_m} = {6 \over d(0)-k} \le {6 \over {k^2\over m}-k} \le
 {12 m\over k^2} < {48m\over N} < 48 N^{-1/2}.
$$
Let
$$
\epsilon = c (mN^{-1/2})^{1/2},
$$
with the positive constant $c$ to be chosen later.
Since $m=o(N^{1/2})$, $\epsilon$ can be made as small as we please and
$$
{3\over N_m} < \epsilon,
$$
if $N$ is large enough.
We also have (since $N_m>{N\over 16m}$)
$$
\epsilon^2 N_m = c^2 {m\over N^{1/2}} N_m >
 {c^2\over 16}\cdot{m\over N^{1/2}}\cdot{N\over m} = {c^2\over 16}N^{1/2},
$$
so that
$$
A \epsilon^2 N_m > A {c^2\over 16}N^{1/2} > k
$$
if $c$ is suitably chosen, i.e., by $Ac^2/32 = 1$.
(Here $A$ is the constant in Lemma \ref{lm:cos_m}.)

Hence the hypothesis of Lemma \ref{lm:cos_m} must fail, and we obtain
(since $N$ is larger than all $\lambda_j$)
$$
N \ge (2-\epsilon)m N_m,
$$
i.e.,
$$
{N\over m} \ge \left(1-c {m^{1/2}\over N^{1/4}}\right) (d(0)-k).
$$
Since $m^{1/2}N^{-1/4}=o(1)$ we have
\begin{eqnarray*}
d(0)-k &\le& \left(1+C{m^{1/2}\over N^{1/4}}\right){N\over m}\\
&\le& \left(1+C{m^{1/2}\over N^{1/4}}\right)
 \left( {k^2\over m} + {2\ell k \over m} + {\ell^2 \over m} \right)\\
&\le& {k^2\over m} + C{k^2 \over m^{1/2}N^{1/4}} + C{\ell k\over m}.\\
\end{eqnarray*}
We also have $k = o({k^2 \over m^{1/2}N^{1/4}})$ since
$m = o(N^{1/2})$ and $k \sim N^{1/2}$.
It follows that
$$
\Abs{\sum_{x\in\ZZ_m}(a(x))^2 - {k^2\over m}} \le
 C\left({k^2\over m^{1/2}N^{1/4}} + {\ell k \over m}\right),
$$
and by Lemma \ref{lm:cst} we obtain (with $k \sim N^{1/2}$)
\begin{eqnarray*}
\Norm{a(x)-{k^2\over m}}_2 &\le&
 C{N^{1/4}\over m^{1/4}} \left( N^{1/4} + {\ell\over m^{1/2}}\right)^{1/2}\\
 &\le&
C \left\{\begin{array}{cl}
  {N^{3/8}\over m^{1/4}} & \mbox{if $\ell\le N^{1/4}m^{1/2}$}\\
  \ & \ \\
  {N^{1/4}\ell^{1/2}\over m^{1/2}} & \mbox{else}
  \end{array} \right.
\end{eqnarray*}
as we had to prove.
\Qed

\mysection{Bibliography}

\end{document}